\documentclass[12pt,leqno]{article}
\usepackage[pdftex]{graphicx}
\usepackage{amsfonts,amssymb,amsmath}
\usepackage{ccaption,mathrsfs,amsmath,amssymb,amsthm}
\usepackage{amsmath, amscd, amsfonts,  tikz, times}

\newcounter{conjecture}
\newcounter{remark}

\newcommand{\eqnsection}{
    \renewcommand{\theequation}{\thesection.\arabic{equation}}
    \makeatletter
    \csname @addtoreset\endcsname{equation}{section}
    \makeatother}
\newtheorem{theorem}{Theorem}

\newtheorem{lemma}{Lemma}

\newtheorem{prop}{Proposition}
\newtheorem{cor}{Corollary}

\newtheorem*{theoremnon}{Theorem}

\newcommand{\dd}{\delta}

\newcommand{\lar}{\longrightarrow}

\newcommand{\aaa}{\alpha}

\newcommand{\CC}{\mathbb{C}}

\def \be{\begin{equation}}
\def \ee{\end{equation}}
\def \bt{\begin{theorem}}
\def \et{\end{theorem}}
\def \bea{\begin{eqnarray}}
\def \eea{\end{eqnarray}}
\def \bas{\begin{eqnarray*}}
\def \eas{\end{eqnarray*}}

\newcommand {\rrr}[1]{(\ref{#1})}

\def \bb{\beta}

\def \ga{\gamma}

\def \si{\sigma}

\def \th{\theta}

\newcommand{\seg}[2]{\stackrel{\line(1,0){#1}}{#2}}
\def \ff{\infty}

\def \AA{{\cal A}}

\def \DD{{\mathbb D}}

\def \GG{{\cal G}}
\def \HH{{\cal H}}
\def \HHH{{\mathbb H}}
\def \II{{\cal I}}

\def \RR{{\mathbb R}}

\def \CCC{{\cal C}}

\def \({\left(}
\def \){\right)}

\def \vski{\vspace{12pt}}

\newcommand{\bsh}{\backslash}

\def \bc{\begin{center} }
\def \ec{\end{center} }
\def \bs{\begin{slide} }
\def \es{\end{slide} }

\def\square{{\vcenter{\vbox{\hrule height.3pt
         \hbox{\vrule width.3pt height5pt \kern5pt
            \vrule width.3pt}
         \hrule height.3pt}}}}
\def\qed{{\hfill $\Box$ \bigskip}}

\newcounter{cccases}
\eqnsection
\begin{document}

\title{The exit time of planar Brownian motion and the Phragm\'en-Lindel\"of principle.}

\author{
\begin{tabular}{c}
\textit{Greg Markowsky} \\
Monash University \\
Victoria 3800, Australia \\
gmarkowsky@gmail.com
\end{tabular}}

\bibliographystyle{amsplain}

\maketitle \eqnsection \setlength{\unitlength}{2mm}

\begin{abstract}
\noindent In this paper a version of the Phragm\'en-Lindel\"of principle is proved using probabilistic techniques. In particular, we will show that if the $p$\textsuperscript{th} moment of the exit time of Brownian motion from a planar domain is finite, then an analytic function on that domain is either bounded by its supremum on the boundary or else goes to $\ff$ along some sequence more rapidly than $e^{|z|^{2p}}$. We also provide a method of constructing domains whose exit time has finite $p$\textsuperscript{th} moment. This allows us to give a general Phragm\'en-Lindel\"of principle for spiral-like and star-like domains, as well as a new proof of a theorem of Hansen. A number of auxiliary results are presented as well.




\end{abstract}

\section{Introduction}

The Phragm\'en-Lindel\"of principle is a method by which the maximum modulus principle can be generalized to certain unbounded domains in $\CC$. The principle roughly states that, on particular domains, analytic functions must either be bounded by their supremum on the boundary of the domain or tend rapidly to $\ff$ along some sequence. We note that, as the principle is generally stated, the precise meaning of "tend rapidly to infinity" will depend upon the domain in question. Our aim in this paper is to prove a general form of the principle using probabilistic arguments, in particular a relationship between the growth of functions and the moments of exit times of planar Brownian motion. We will also show how the principle can be applied in a number of special cases.\\

In order to give a precise statement of our main result, we need a few definitions. In what follows, $B_t$ will always refer to a planar Brownian motion. For any domain $W \subseteq \CC$ we let $T_W = \inf\{t \geq 0: B_t \notin W\}$ be the first exit time of Brownian motion from $W$. The notations $E_a$ and $P_a$ will be used to refer to expectation and probability conditioned upon $B_0 = a$ a.s. If $E_a[T_W^p]< \ff$ for some $a \in W, p>0$, then the connectedness of $W$ implies that $E_b[T_W^p]< \ff$ for all $b \in W$ (see \cite[(3.13)]{burk}), and we will in this case simply write $E[T_W^p] < \ff$. $\dd W$ denotes the boundary of $W$ in $\CC$; that is, $\dd W$ does not include the point at $\ff$. We will prove the following theorem.

\begin{theorem} \label{phraglind1}
Let $W$ be a domain such that $E[T^p_W] < \ff$. Suppose that $f$ is an analytic function on $W$ such that $\limsup_{z \lar \dd W}|f(z)| \leq K < \ff$, and $|f(z)| \leq Ce^{C|z|^{2p}} + C$ for some $C>0$. Then $|f(z)| \leq K$ for all $z \in W$.
\end{theorem}

The proof will be given in Section \ref{p-l}. The theorem encompasses some well-known special cases, as well as some which appear to be new, as will be shown in Section \ref{app}. For instance, formulations are available for an infinite wedge and arbitrary simply connected domains, as well as for general star-like and spiral-like domains. Further formulations are possible which make use of a method, presented in Section \ref{starspi}, of building domains whose exit time has finite $p$\textsuperscript{th} moment.


\section{Proof of Theorem \ref{phraglind1}.} \label{p-l}

The key to our investigation will be the following pair of results, which we will collectively refer to as {\it Burkholder's theorem}.

\begin{theoremnon}[{\bf Burkholder}]
(i) For any $p \in (0,\ff)$ there are constants $c_p, C_p > 0$ such that for any stopping time $\tau$ we have

\begin{equation} \label{bteq}
c_pE_a[(\tau + |a|^2)^p] \leq E_a[|B^*_\tau|^{2p}] \leq C_p E_a[(\tau + |a|^2)^p].
\end{equation}

In particular, $E_a[\tau^p] < \ff$ if, and only if, $E_a[|B^*_\tau|^{2p}]< \ff$. \\

(ii) For any $p \in (0,\ff)$ there is a constant $C_p>0$ such that for any stopping time $\tau$ with $E_a[\ln \tau]< \ff$ we have

\begin{equation} \label{}
E_a[|B_\tau|^{2p}] \leq E_a[|B^*_\tau|^{2p}] \leq C_p E_a[|B_\tau|^{2p}].
\end{equation}

\end{theoremnon}

It may be tempting to see part $(i)$ at least as a straightforward consequence of the standard Burkholder-Davis-Gundy inequality by separately bounding the supremums of the real and imaginary parts of $B_t$; however, this argument is not quite valid, since a stopping time for $B_t$ need not be a stopping time for its projection onto the real or imaginary axis. The reader who would like to see a proof of the theorem is therefore referred to \cite{burk}. \\

Before proving Theorem \ref{phraglind1}, we give a preliminary result on subharmonic functions. In what follows, $cl(W)$ will denote the closure of the set $W$ in $\CC$.

\begin{prop} \label{calig}
Let $W$ be a domain with $E[T_W^p] < \ff$. Suppose that $u$ is a continuous function on $cl(W)$ which is subharmonic on $W$ and satisfies $\sup_{z \in \dd W} u(z) \leq K$, for some $K > 0$. Suppose further that $u(z) \leq C|z|^{2p} + C$ for some $C < \ff$. Then $u(z) \leq K$ for all $z \in W$.
\end{prop}

{\bf Proof:} Let $S_M = \inf\{t \geq 0:|B_t|=M\}$, and fix $a \in W$. Since $u$ is subharmonic, $u(a) \leq E_a[u(B_{T_W \wedge S_M})]$ (see \cite[Sec.'s 2.IV.3 and 2.IX.3]{doobook}). We would like to let $M \lar \ff$ to obtain $u(a) \leq E_a[u(B_{T_W})]\leq K$; note that the conditions on $u$ imply that $|u(B_{T_W \wedge S_M})| \leq C|B_{T_W \wedge S_M}|^{2p} + C \leq C|B^*_{T_W}|^{2p} + C$, and $E_a[C|B^*_{T_W}|^{2p} + C] < \ff$ by Burkholder's Theorem. The dominated convergence theorem therefore applies, and we get $u(a) \leq \lim_{M \nearrow \ff} E_a[u(B_{T_W \wedge S_M})] = E_a[u(B_{T_W})] \leq K$. \qed


{\bf Proof of Theorem \ref{phraglind1}:} We may assume $K=1$. Set

\begin{equation} \label{17dd}
\log^+ x = \begin{cases} \log x & x > 1,\\
0 & x  \leq 1.
\end{cases}
\end{equation}

The function $u(z) = \log^+ |f(z)|$ is the maximum of two subharmonic function, and is therefore subharmonic. Note that the conditions on $f$ imply that $\sup_{z \in \dd W} u(z) = 0$ and $u(z) \leq C|z|^{2p} + C$ for some (possibly different) $C>0$. Applying Proposition \ref{calig} now implies that $u(z) \leq 0$ for all $z \in W$, and the result follows. \qed




\section{Applications} \label{app}

In order to state useful special cases of Theorem \ref{phraglind1}, we need to find domains $W$ for which $E[T^p_W] < \ff$. The primary method for doing this was also exhibited by Burkholder in \cite{burk}, and involves the {\it Hardy norm} of conformal maps, as we now describe. If $W \subsetneq \CC$ is simply connected, then the Riemann Mapping Theorem guarantees the existence of a conformal map $f_a$ mapping $\DD$ onto $W$ which takes $0$ to $a$. The Hardy norm $|| \cdot ||_{H^{2p}}$ of $f_a$ is defined as

\begin{equation} \label{h2def}
||f_a||_{H^{2p}} := \Big(\sup_{r<1} \frac{1}{2\pi} \int_{0}^{2\pi} |f_a(re^{i \th})|^{2p} d\th \Big)^{1/2p} .
\end{equation}

The map $f_a$ is not uniquely determined; however any two such maps differ only by precomposition with a rotation, so the value of $||f_a||_{H^{2p}}$ is independent of the choice of $f_a$. In light of this observation, let us set $\HH^a_{2p}(W) = ||f_a||_{H^{2p}}$, and note that $\HH^a_{2p}$ provides a sort of measure on the size and shape of domains. The following result also first appeared in \cite{burk}.

\begin{prop} \label{hardyequiv}
For any $p \in (0,\ff)$ there are constants $c_p, C_p > 0$ such that if $W \subsetneq \CC$ is simply connected then

\begin{equation} \label{bteq2}
c_pE_a[(T_W + |a|^2)^p] \leq \HH^a_{2p}(W)^{2p} \leq C_p E_a[(T_W + |a|^2)^p].
\end{equation}

In particular, $E[T_W^p] < \ff$ if, and only if, $\HH^a_{2p}(W)< \ff$ for any $a \in W$.
\end{prop}

To obtain our first variant of Theorem \ref{phraglind1}, we remark that it is known that if $W \subsetneq \CC$ is simply connected then $E[T^p_W] < \ff$ for any $p < \frac{1}{4}$; this is proved in \cite[p. 301]{burk}, and the fact that the $H^{2p}$-Hardy norm of the Koebe function $f(z)=\frac{z}{(1-z)^2}$, which maps $\DD$ conformally onto $\CC \bsh (-\ff, 1/4]$, is finite if and only if $p<1/4$ shows it cannot be improved (this is also implied by Theorem \ref{lamoreau} below). We obtain

\begin{cor} \label{simpcon}
Suppose that $f$ is an analytic function on a simply connected domain $W \subsetneq \CC$ such that $\limsup_{z \lar \dd W}|f(z)| \leq K < \ff$, and $|f(z)| \leq Ce^{C|z|^{2p}} + C$ for some $p < \frac{1}{4}$. Then $|f(z)| \leq K$ for all $z \in W$.
\end{cor}

A domain $W$ is {\it spiral-like of order $\si \geq 0$ with center $a$} if, for any $z \in W$, the spiral $\{a+(z-a) \mbox{ exp}(te^{-i \si}) : t \leq 0\}$ also lies within $W$ (cf. \cite{space}). A natural question in light of Theorem \ref{phraglind1} can be posed: for a given spiral-like domain $W$ and $p>0$, is $E[T_W^p]<\ff$? In \cite{hansenspi}, Hansen gave a geometric condition for the finiteness of $\HH^a_{2p}(W)$, but before stating the result let us examine the question in more detail. There is no loss of generality in assuming $a=0$, and we will do so henceforth. Hansen showed that the key quantity for our purposes is the measure of the largest arc in the set $W \cap \{|z|=r\}$ (taken as a set on the circle), and with this in mind we let

\begin{equation} \label{bigmax}
\AA_{r,W} = \max \{m(E): E \mbox{ is a subarc of } W \cap \{|z|=r\}\},
\end{equation}

where $m$ denotes angular Lebesgue measure on the circle. Spiral-likeness implies that $\AA_{r,W}$ is nondecreasing in $r$, so we may let $\AA_W = \lim_{r \nearrow \ff} \AA_{r.W}$. Hansen's result is as follows, translated from an analytic statement into the corresponding probabilistic one via Proposition \ref{hardyequiv}.

\begin{theorem} \label{lamoreau}
If $W$ is a spiral-like domain of order $\si$ with center 0, then $E[T_W^p] < \ff$ if, and only if, $p < \frac{\pi}{2 \AA_W \cos^2 \si}$.
\end{theorem}

\vspace{.02in}

\begin{cor} \label{page}
Suppose $W$ is a spiral-like domain of order $\si$ with center $0$. If $f$ is an analytic function on $W$ such that $\limsup_{z \lar \dd W}|f(z)| \leq K < \ff$, and $|f(z)| \leq Ce^{C|z|^{2p}} + C$ for any $p < \frac{\pi}{2 \AA_W \cos^2 \si}$, then $|f(z)| \leq K$ for all $z \in W$.
\end{cor}

Star-like domains have been studied more intensively by analysts than spiral-like ones. A domain $W$ is called {\it star-like with center $a$} if the line segment connecting $a$ to $z$ lies within $W$ for every $z \in W$. Note that a star-like domain is simply spiral-like of order $\si =0$, and Corollary \ref{page} therefore takes the following form as a special case.

\begin{cor} \label{kristen}
Suppose $W$ is a star-like domain with center $0$. If $f$ is an analytic function on $W$ such that $\limsup_{z \lar \dd W}|f(z)| \leq K < \ff$, and $|f(z)| \leq Ce^{C|z|^{2p}} + C$ for any $p < \frac{\pi}{2 \AA_W}$, then $|f(z)| \leq K$ for all $z \in W$.
\end{cor}

Note that convex domains are trivially star-like, and the previous theorem therefore applies to any convex domains as well. Let us now set $N_\aaa = \{re^{i\th}: r \in (0,\ff), \th \in (\frac{-\aaa}{2},\frac{\aaa}{2})\}$; $N_\aaa$ is the angular wedge with vertex at 0 which is symmetric about the real axis and has angular width $\aaa$. $N_\aaa$ is star-like, and Corollary \ref{kristen} therefore reduces further to the following.

\begin{cor} \label{edgewedge}
Suppose that $f$ is an analytic function on $N_\aaa$ such that $\limsup_{z \lar \dd W}|f(z)| \leq K < \ff$, and $|f(z)| \leq Ce^{C|z|^{2p}} + C$ for some $p < \frac{\pi}{2\aaa}$. Then $|f(z)| \leq K$ for all $z \in N_\aaa$.
\end{cor}

This is a commonly stated form of the Phragm\'en-Lindel\"of principle (see \cite[Sec. VI.4]{conway}, for instance), and was in fact previously given a different proof using Brownian motion in \cite[Thm. V.1.8 and Ex. V.2]{bassy}. It may also be noted that this example shows that Corollary \ref{page} (and, in turn, Theorem \ref{phraglind1}) is sharp in the sense that the function $f(z) = e^{z^{\pi/\aaa}}$ is analytic on $N_\aaa$ and bounded in modulus by 1 on $\dd N_\aaa$, but is clearly unbounded on $N_\aaa$.\\

It should be mentioned that, due to the nature of the random variable $T_W$, it is trivial that if $W_1 \subseteq W_2$ then $T_{W_1} \leq T_{W_2}$ a.s. The conclusions of Corollaries \ref{simpcon} through \ref{edgewedge} therefore hold if the domains in question are replaced by smaller ones, which need not be simply connected.

\section{A method for generating domains whose exit time have finite $p$\textsuperscript{th} moments and a probabilistic proof of Theorem \ref{lamoreau}} \label{starspi}

It is clear that in order to find other applications of Theorem \ref{phraglind1} one must be able effectively bound moments of Brownian exit times. The natural initial attempt in this direction might be to reduce complicated domains in some way to simpler ones for which we have good bounds; in particular we might hope that if $V, W$ are domains with $E[T_V^p]< \ff$ and $E[T_W^p]< \ff$ then we can construct a new domain out of the two of them whose exit times has also a finite $p$-th moment. However, it is immediately clear that we may not simply take the union of any such $V$ and $W$, since for instance if $V = \{Re(z) >0\}$ and $W = \{Re(z) <1\}$, then $E[T_V^p], E[T_W^p] < \ff$ for $p < 1/2$, but $V \cup W = \CC$, whose exit time is infinite.

\vski

In this section we describe a method of circumventing this difficulty in order to build domains whose exit times have a finite $p$-th moment. 
Given two domains, $V$ and $W$, we let $\dd V^+ = \dd V \cap W$ and $\dd W^+ = \dd W \cap V$, as is shown. 

\hspace{1.4in} \includegraphics[height=3in]{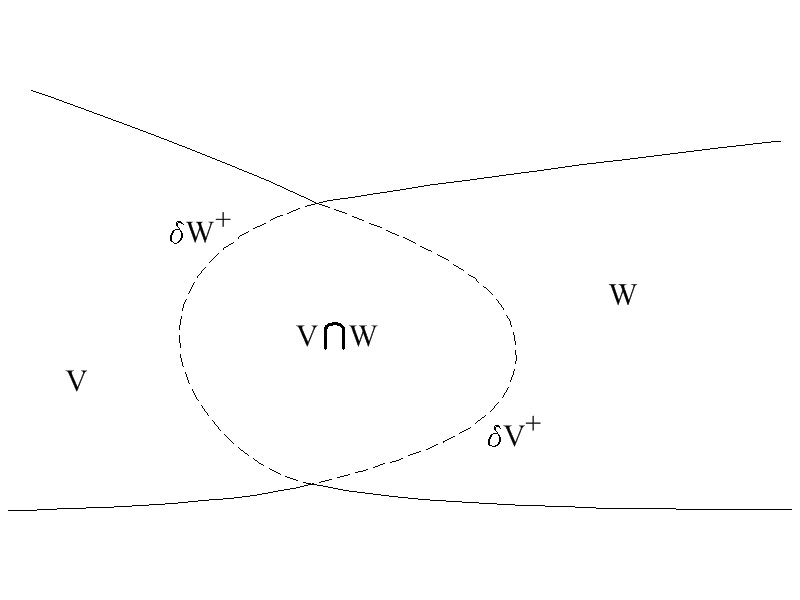}

Note that $\dd(V \cup W) = (\dd V \cup \dd W) \bsh (\dd V^+ \cup \dd W^+)$, so that a Brownian motion exits $V \cup W$ precisely when it hits a boundary point of $V$ or $W$ which is not contained in $\dd V^+$ or $\dd W^+$. The following theorem and the lemmas which follow constitute the aforementioned technique for building domains while keeping control of the moments of exit time.

\begin{theorem} \label{flydag}
Suppose that $V$ and $W$ are domains with nonempty intersection, neither of which is contained in the other. Suppose further that $E[T_V^p]< \ff$ and $E[T_W^p]< \ff$. Let $\dd V^+$ and $\dd W^+$ be defined as above, and assume that the following conditions are satisfied:

\begin{itemize} \label{}


\item[(i)] $\sup_{a \in \dd V^+} E_a[T^p_W] < \ff$;

\item[(ii)] $\sup_{a \in \dd W^+} E_a[T^p_V] < \ff$;

\item[(iii)] $\sup_{a \in \dd V^+} P_a(B_{T_W} \in \dd W^+) < 1$.

\end{itemize}

Then $E[T_{V \cup W}^p]< \ff$.
\end{theorem}

{\bf Proof:} The conditions on $V$ and $W$ imply that $\dd V^+$ and $\dd W^+$ are both nonempty. Form a sequence of stopping times $\tau_n,\tau'_n$ as follows. Let $\tau_0 = 0$, $\tau'_1 = \inf \{t \geq \tau_0:B_t \in W^c\}$, and $\tau_1 = \inf \{t \geq \tau'_1:B_t \in V^c\}$. Continue in this manner, letting $\tau'_n = \inf \{t \geq \tau_{m-1}:B_t \in W^c\}$ and $\tau_n = \inf \{t \geq \tau'_n:B_t \in V^c\}$. Note that $\tau_n \leq T_{V \cup W}$ for all $n$, since if $\tau_n = T_{V \cup W}$ or $\tau'_n = T_{V \cup W}$ for some $n$ then $\tau_m=\tau'_m = T_{V \cup W}$ for all $m > n$. Furthermore, if we set $\tau_\ff=\lim_{n \lar \ff} \tau_n=\lim_{n \lar \ff} \tau'_n$, then on the event $\{ \tau_\ff < \ff\}$ by continuity we will have $B_{\tau_\ff} \in V^c$ and $B_{\tau_\ff} \in W^c$, so that $T_{V \cup W} = \tau_\ff$. It follows from this that $\tau_\ff = T_{V \cup W}$ a.s. We will start a Brownian motion at a point $a \in \dd V^+$, and write

\begin{equation} \label{portland}
\tau_\ff = \sum_{n=1}^{\ff} (\tau'_n - \tau_{n-1}) + \sum_{n=1}^{\ff} (\tau_n - \tau'_n).
\end{equation}

Note that the strong Markov property implies that

\begin{equation} \label{grape}
\begin{gathered}
E_a[(\tau'_n - \tau_{n-1})^p] \leq P_a(\tau_{n-1}<T_{V \cup W})\sup_{a \in \dd V^+} E_a[T^p_W], \\ E_a[(\tau_n - \tau'_n)^p] \leq P_a(\tau'_n < T_{V \cup W})\sup_{a \in \dd W^+} E_a[T^p_V].
\end{gathered}
\end{equation}

Now, $\{\tau'_n<T_{V \cup W}\} \subseteq \{\tau_{n-1}<T_{V \cup W}\} \subseteq \{\tau'_{n-1}<T_{V \cup W}\}$, yielding monotonicity in the corresponding probabilities; but in fact more is true, namely that $P_a(\tau'_n<T_{V \cup W}) \leq \sup_{a \in \dd V^+} P_a(B_{T_W} \in \dd W^+) P_a(\tau_{n-1}<T_{V \cup W})$, again by the strong Markov property and the fact that $\tau'_n<T_{V \cup W}$ precisely when $B_{\tau'_n} \in \dd W^+$. It follows that, with $r = \sup_{a \in \dd V^+} P_a(B_{T_W} \in \dd W^+)$, we have $P_a(\tau'_n<T_{V \cup W}) \leq r^n P_a(\tau_0 < T_{V \cup W}) = r^n$, and we have also $P_a(\tau_{n-1}<T_{V \cup W}) \leq P_a(\tau'_{n-1}<T_{V \cup W}) \leq r^{n-1}$. Employing this bound, \rrr{grape}, and Minkowski's inequality, \rrr{portland} becomes

\begin{equation} \label{}
\begin{split}
||\tau_\ff||_p & \leq \sum_{n=1}^{\ff} ||\tau'_n - \tau_{n-1}||_p + \sum_{n=1}^{\ff} ||\tau_n - \tau'_n||_p \\
& \leq \sum_{n=1}^{\ff} r^{(n-1)/p} \sup_{a \in \dd V^+} E_a[T^p_W]^{1/p} + \sum_{n=1}^{\ff} r^{n/p} \sup_{a \in \dd W^+} E_a[T^p_V]^{1/p}.
\end{split}
\end{equation}

Assumptions $(i)-(iii)$ show this quantity to be finite. \qed

It may appear that the conditions in Theorem \ref{flydag} are difficult to check, but in fact they are quite easy to check in many cases where the domains in question have particularly nice boundaries; the ensuing two lemmas provide simple checks which are sufficient for our purposes. A {\it Jordan domain} is a bounded, simply connected domain whose boundary is homeomorphic to a circle. We will call a domain $W$ a {\it Jordan$^*$ domain} if there is a Jordan domain $U$ and and M\"obius transformation $\phi$ such that $\phi(U)=W$; this essentially gives us a class of domains with the same nice properties as Jordan domains but which now include many unbounded domains. The following lemma (which, naturally, holds with the roles of $V$ and $W$ interchanged), allow us to confirm conditions $(i)$ and $(ii)$ in many instances.


\begin{lemma} \label{ziyi}
Suppose $W$ is a Jordan$^*$ domain with $E[T_{W}^p]<\ff$, and $\dd V^+$ is bounded. Then $\sup_{a \in \dd V^+} E_a[T^p_W] < \ff$.
\end{lemma}

{\bf Proof:} Since $W$ is simply connected, Proposition \ref{hardyequiv} allows us to consider the quantity $\HH^a_{2p}(W)$ in place of $E_a[T^p_W]$. Let $\GG(a) = \HH^a_{2p}(W)$ for $a \in W$ and $\GG(a) = |a|$ for $a \in \dd W$. We will show that $\GG$ is continuous on the closure $cl(W)=W \cup \dd W$. Suppose that $a =: a_0 \in W$, and that $\{a_n\}_{n=1}^\ff$ is a sequence of points in $W$ converging to $a$. Let $f_{a_0}$ be defined as for \rrr{h2def}; since $W$ is a Jordan$^*$ domain, Caroth\'eodory's Theorem (\cite[Thm. 5.1.1]{krantz}) implies that $f$ extends to a homeomorphism (taking values in the sphere $\CC \cup \{\ff\}$) from $\dd \DD$ to $\dd W$, and our assumptions on $W$ show that $f$, so extended, is in $L^{2p}(\dd \DD)$. Let $b_n = f^{-1}(a_n)$ for each $n$, and let $\phi_n(z) = \frac{z+b_n}{1 + \seg{1.2}{b_n} z}$. $\phi_n$ is the disk automorphism taking $0$ to $b_n$, so that $f_{a_n} := f_{a_0} \circ \phi_n$ is a conformal map from $\DD$ onto $W$ taking $0$ to $a_n$. We have

\begin{equation} \label{}
\begin{split}
\GG(a_n) = ||f_{a_n}||_{H^{2p}} & = \Big(\frac{1}{2\pi} \int_{0}^{2\pi} |f_{a_n}(e^{i \th})|^{2p} d\th \Big)^{1/2p} \\
& = \Big(\frac{1}{2\pi} \int_{0}^{2\pi} |f_{a_0}(\phi_n(e^{i \th}))|^{2p} d\th \Big)^{1/2p} \\
& = \Big(\frac{1}{2\pi} \int_{0}^{2\pi} |f_{a_0}(e^{i \th})|^{2p} |\phi_n'(e^{i \th})| d\th \Big)^{1/2p} \\
& = \Big(\frac{1}{2\pi} \int_{0}^{2\pi} |f_{a_0}(e^{i \th})|^{2p} \frac{1-|b_n|^2}{|1 + \seg{1.2}{b_n} z|^2} d\th \Big)^{1/2p}.
\end{split}
\end{equation}

As $n \lar \ff$, $a_n \lar a_0$, which means that $b_n \lar 0$ and the kernels $K_n = \frac{1}{2\pi}\frac{1-|b_n|^2}{|1 + \seg{1.2}{b_n} z|^2}$ approach $\frac{1}{2\pi}$ uniformly. It follows easily from this that $\GG(a_n) \lar \GG(a_0) = \GG(a)$, so that $\GG$ is continuous at $a$. Now suppose that $a \in \dd W$, and that $\{a_n\}_{n=1}^\ff$ is again a sequence of points in $W$ converging to $a$. Choose $a_0$ be chosen arbitrarily in $W$, and let the $b_n$'s and $\phi_n$'s be defined as before. We will now have $b_n \lar b = f_{a_0}^{-1}(a) \in \dd \DD$ as $n \lar \ff$. The kernels $K_n$ are positive, have total mass 1, and approach $0$ uniformly on $\dd \DD \backslash \{b\}$; they therefore form a Dirac sequence (see \cite[Ch. XI]{lang}), and it follows by standard methods that

\begin{equation} \label{}
\begin{split}
\GG(a_n) = ||f_{a_n}||_{H^{2p}} \lar |f_{a_0}(b)| = |a| = \GG(a).
\end{split}
\end{equation}

$\GG$ is therefore continuous on all of $cl(W)$. As such, it must remain bounded on any compact set, and the result follows. \qed


The following lemma is useful for checking condition $(iii)$ in Theorem \ref{flydag} in certain instances. Recall that a curve $\ga$ is {\it analytic} if at every point $v$ on $\ga$ there is a neighborhood $U$ of $v$ and a conformal map $\phi$ from $\DD$ onto $U$ such that $\phi(\DD \cap \RR) = \ga \cap U$.

\begin{lemma} \label{page}
Suppose that $W$ is a Jordan* domain, with two points $v_1, v_2 \in \dd W$ which both lie in boundary arcs which are analytic. Suppose that $\CCC$ is a simple curve lying in $W$ which connects $v_1, v_2$; then $W \backslash \CCC$ has two components, $W_1$ and $W_2$. Suppose further that $\CCC$ is differentiable at $v_1, v_2$, and the angles $\CCC$ makes with the boundary arcs at $v_1,v_2$ are not zero. Then $\sup_{a \in \CCC} P_a(B_{T_W} \in \dd W_1) < 1$, and likewise $\sup_{a \in \CCC} P_a(B_{T_W} \in \dd W_2) < 1$.
\end{lemma}

\hspace{1.4in} \includegraphics[height=3in]{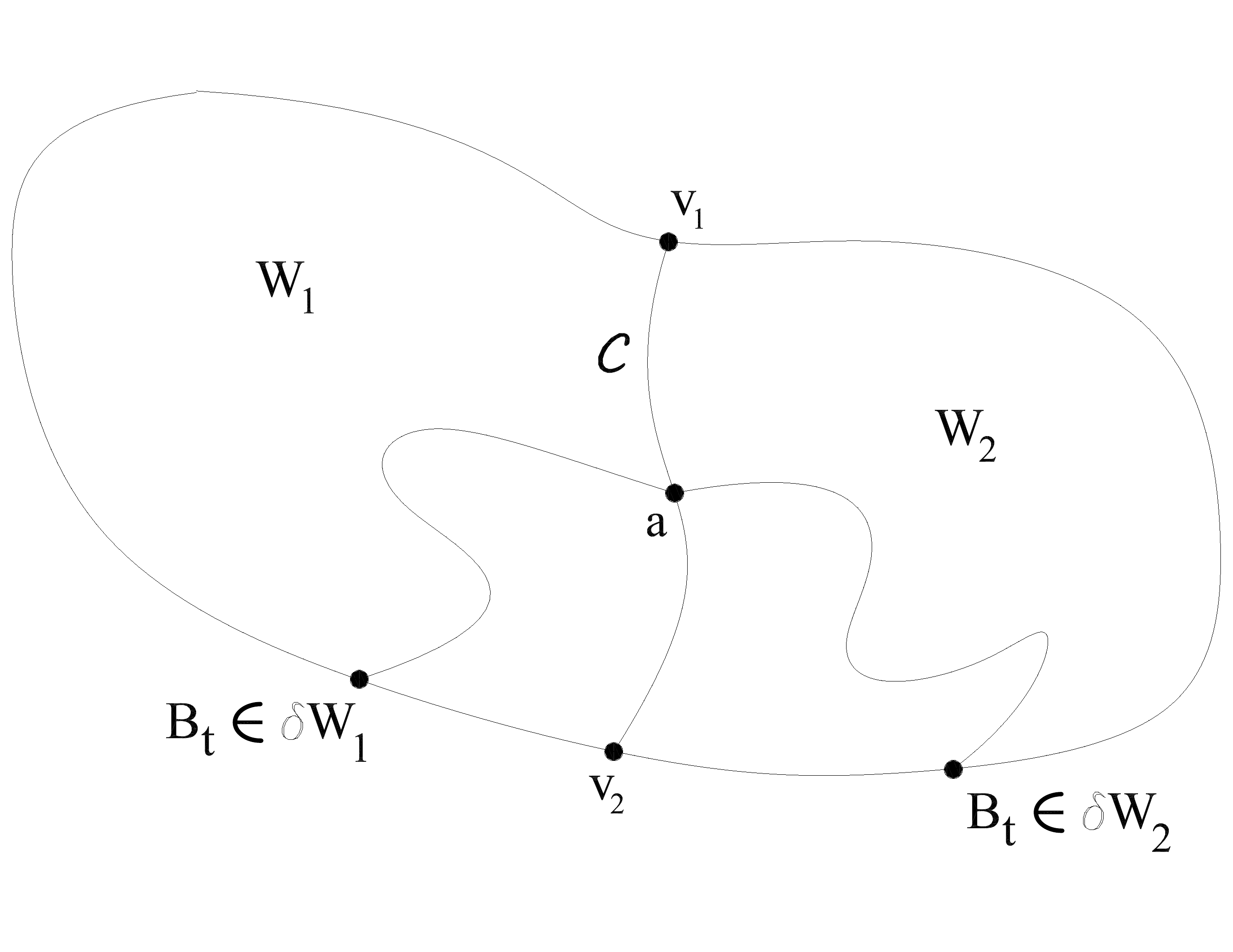}

{\bf Proof:} That $W \backslash \CCC$ consists of two components is a simple consequence of the Jordan Curve Theorem. Let $\phi$ be a conformal map from $W$ to $\HHH = \{Im(z)>0\}$ taking $v_1, v_2$ to $0, \ff$. Then Carath\'eodory's Theorem assures us of a continuous extension of $\phi$ to the boundary, mapping $\dd W$ to $\RR$. Replace $\phi$ be $-1/\phi$ if necessary so that $\phi$ takes $\dd W_2 \cap \dd W$ to $\RR^+ = \{Im(z)=0, Re(z)>0\}$, and $\dd W_1 \cap \dd W$ to $\RR^- = \{Im(z)=0, Re(z)<0\}$. The conformal invariance of Brownian motion shows that $\sup_{a \in \CCC} P_a(B_{T_W} \in \dd W_1) = \sup_{a \in \phi(\CCC)} P_a(B_{T_\HHH} \in \RR^-)$, and this is equal to $\sup_{a \in \phi(\CCC)} \frac{Arg(a)}{\pi}$, with $Arg$ denoting the principle branch of the argument function; this is a simple consequence of the fact that $a \lar \frac{Arg(a)}{\pi}$ is harmonic, bounded, and equal to $1$ on $\RR^-$ and $0$ on $\RR^+$, and is therefore equal to the harmonic measure of $\RR^-$ on $\HHH$. Now, $\phi(\CCC)$ is a curve connecting $0$ to $\ff$, and the Schwarz reflection principle (see \cite[Sec. I.1.6]{dur}) shows that $\phi$ extends to be analytic with nonzero derivative at $v_1, v_2$, so that $\phi(\CCC)$ meets $\dd \HHH$ at nonzero angles at $0, \ff$. This implies that if $a$ approaches $0$ or $\ff$ along $\phi(\CCC)$, $Arg(a)$ remains bounded away from 1, and it follows by a compactness argument that $\sup_{a \in \phi(\CCC)} \frac{Arg(a)}{\pi}< 1$. This shows that $\sup_{a \in \CCC} P_a(B_{T_W} \in \dd W_1) < 1$, and the corresponding statement for $\dd W_2$ follows upon interchanging the roles of $W_1$ and $W_2$. \qed

Let us now prove Hansen's Theorem (Theorem \ref{lamoreau}). Suppose $W$ is spiral-like of order $\si \geq 0$ with center $0$. If a point $z = re^{i \aaa}$ lies in $W^c$, then in fact the entire curve $\{re^{i\aaa} \mbox{ exp}(te^{-i \si}) : t \geq 0\}$ must lie in $W^c$ as well, and it therefore suffices to consider domains of the form $S_{r,D} = \CC \bsh \{re^{i\aaa} \mbox{ exp}(te^{-i \si}): t \geq 0, \aaa \in D\}$, where $r>0$ and $D$ is a finite subset of $[0,2\pi)$, and then to approximate arbitrary spiral-like domains by domains of this form. Furthermore, standard Brownian scaling and rotation invariance allows us to assume $0 \in D$ and to consider only $S_D := S_{1,D}$. Let $D$ be given as $\{\aaa_0, \aaa_1, \ldots , \aaa_k\}$, with $0 = \aaa_0 < \aaa_1 < \aaa_2 < \ldots < \aaa_k < \aaa_{k+1} := 2\pi $, and let us extend our prior notation for the infinite wedge by defining $N_{\aaa}^{\bb} = \{e^{i\th} \mbox{ exp}(te^{-i \si}) : t \in \RR, \th \in (\aaa,\bb)\}$ for $\aaa < \bb$; we may therefore write $S_D = \DD \cup \Big(\cup_{j=0}^k N_{\aaa_j}^{\aaa_{j+1}}\Big)$. The following picture gives an example of an $S_D$ divided up into such a union, where for simplicity we have set $\si=0$, corresponding to the star-like case.

\hspace{1.4in} \includegraphics[height=4in,width=4in]{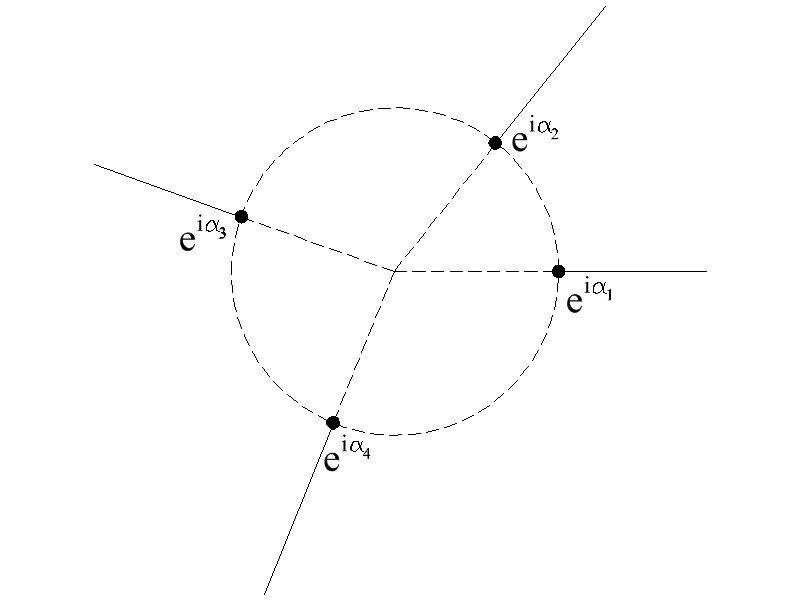}

Using the notation from Section \ref{app}, we have $\AA_{S_D} = \max_j (\aaa_{j+1}-\aaa_j)$. It is shown in \cite[p. 280]{hansenspi} that the function

\begin{equation} \label{}
f(z) = e^{i (\aaa_{j+1}-\aaa_j)/2} \mbox{exp} \Big[\frac{(\aaa_{j+1}-\aaa_j) \cos \si}{\pi} e^{-i(\aaa_{j+1}-\aaa_j)}\mbox{Log} \Big(\frac{1+z}{1-z}\Big)\Big]
\end{equation}


maps $\DD$ conformally onto $N_{\aaa_j}^{\aaa_{j+1}}$, and that $||f||_{H^2p} < \ff$ precisely when $p < \frac{\pi}{2 (\aaa_{j+1}-\aaa_j) \cos^2 \si}$. Thus, if $p \geq \frac{\pi}{2 \AA_{S_D} \cos^2 \si}$, then Proposition \ref{hardyequiv} implies that $E[T_{N_{\aaa_j}^{\aaa_{j+1}}}]=\ff$ for some $j$, and by monotonicity we have $E[T_{S_D}]=\ff$. On the other hand, if $p < \frac{\pi}{2 \AA_{S_D} \cos^2 \si}$, then $E[T_{N_{\aaa_j}^{\aaa_{j+1}}}]<\ff$ for all $j$, and we may construct $S_D$ by adding the domains $N_{\aaa_j}^{\aaa_{j+1}}$ one by one to the disk $\DD$. Since $\DD$ is bounded and $\dd \DD$ intersects the boundary of each $N_{\aaa_j}^{\aaa_{j+1}}$ along an analytic boundary arc and with nonzero angle, Theorem \ref{flydag} can be applied via Lemmas \ref{ziyi} and \ref{page} to conclude that the $p$-th moment is bounded at each step, and therefore for the full domain $S_D$. This completes the proof of Theorem \ref{lamoreau}. \\

Naturally, the method outlined in this section can be applied in many instances in order to form a domain which is not spiral-like, or even simply connected, and whose exit time has finite $p$\textsuperscript{th} moment.

\section{Acknowledgements}

I'd like to thank Andrea Collevecchio and Burgess Davis for helpful conversations. I am also grateful for support from Australian Research Council Grants DP0988483 and DE140101201.

\bibliographystyle{alpha}
\bibliography{CABMbib}

\end{document}